\documentclass{elsarticle}

\usepackage{amsmath}
\usepackage{bm}

\usepackage{tikz}
\usepackage{pgfplots}
\usepgflibrary{arrows}
\usepackage{rotating} 
\graphicspath{{pictures/}}

\newtheorem{thm}{Theorem}

\newproof{pf}{Proof}

\journal{arXiv.org} 

\begin{document}

\begin{frontmatter}

\title{Computational identification of the lowest space-wise dependent coefficient of a parabolic equation}

\author[nsi,sffu]{Petr~N.~Vabishchevich\corref{cor}}
\ead{vab@ibrae.ac.ru}

\address[nsi]{Nuclear Safety Institute, Russian Academy of Sciences, 52, B.~Tulskaya, Moscow, Russia}

\address[sffu]{North-Eastern Federal University, 58, Belinskogo, Yakutsk, Russia}

\cortext[cor]{Corresponding author}

\begin{abstract}
In the present work, we consider a nonlinear inverse problem of identifying the lowest coefficient of a parabolic equation. 
The desired coefficient depends on spatial variables only. Additional information about the solution is given 
at the final time moment, i.e., we consider the final redefinition.
An iterative process is used to evaluate the lowest coefficient, where at each iteration 
we solve the standard initial-boundary value problem for the parabolic equation.
On the basis of the maximum principle for the solution of the differential problem,
the monotonicity of the iterative process is established along with the fact that
the coefficient approaches from above. 
The possibilities of the proposed computational algorithm are illustrated by numerical examples 
for a model two-dimensional problem.
\end{abstract}


\end{frontmatter}

\section{Introduction}

Mathematical modeling of many applied problems of science and engineering
leads to numerical solving inverse problems for equations with
partial derivatives \cite{alifanov2011inverse,lavrentev1986ill}.
In the theoretical study of such problems, the main attention is given to
issues of well-posedness of problems, the uniqueness of the solution and its stability.

For parabolic equations, inverse coefficient problems attract particular interest.
In these problems, identification of coefficients of equations and/or their right-hand side 
is conducted using some additional information about the solution. It is possible to identify dependence
of coefficients on time or on spatial variables~\cite{isakov2017inverse,prilepko2000methods}.
Problems of identifying the right-hand side of the equation belong to the class of linear inverse problems. 
Other inverse coefficient problems are non-linear that complicate significantly their study.

Among inverse problems of coefficient identification for parabolic equations we can highlight
problems of determining the dependence of the lowest coefficient (reaction coefficient) on spatial variables. 
As a rule, additional conditions are formulated as the solution value at the final time moment and so, in this case, 
we speak of the final redefinition.
In a more general case, a redefinition condition is formulated as some integral time-average relation 
(integral redefinition). The existence and uniqueness of the solution to such an inverse problem and 
the well-posedness of this problem are considered in a number of works. 
The pioneer works~\cite{isakov1991inverse,prilepko1987solvability}
are devoted to problems with the final redefinition in H\"{o}lder classes and are based on the Schauder principle.
Later works~\cite{kamynin2013inverse,kamynin2010two,prilepko1992inverse,prilepko1993inverse} deal
with problems with integral redefinition and so, they are studied in Sobolev classes.

In works~\cite{kostin2015inverse,prilepko1992inverse} 
(see also~\cite[Theorem 9.1.4]{isakov2017inverse}),
the existence of the solution to the inverse problem of finding the lowest coefficient of a parabolic equation  
is proved constructively. Namely, an iterative process is used with solving the standard initial-boundary 
parabolic problem at each iteration. It seems natural to implement this approach 
in a corresponding computational algorithm.

The standard approach to numerical solving inverse coefficient problems for partial differential equations 
is associated with the minimization of the residual functional using  regularization 
procedures~\cite{samarskii2007numerical,vogel2002computational}.
Computational algorithms are based on the employment of gradient iterative methods, where
we solve initial-boundary value problems both for the initial parabolic equation
and the equation that is conjugate to it. For problems of identifying the lowest coefficient of parabolic equations, 
which depends only on spatial variables, the optimization method
in a combination with finite element approximations in space is used in the work~\cite{chen2006solving}. 
Among later studies in this direction, we mention \cite{cao2018reconstruction,yang2008inverse}.

In the work~\cite{liu2011lie}, an iterative process for the identification of the reaction coefficient 
in the diffusion-reaction equation is proposed without any exact mathematical justification. 
For model one- and two-dimensional boundary value problems, using finite-difference approximations in space, 
the efficiency of this computational algorithm has been demonstrated.
This approach has been also applied to some other inverse problems for parabolic equations,
in particular, for identifying the highest coefficient~\cite{liu2008lgsm}. 

In the present paper, we construct a computational algorithm for identifying the lowest coefficient 
with the final redefinition, which is based on an iterative adjustment of the reaction coefficient 
similarly to~\cite{isakov2017inverse,kostin2015inverse}.
The main attention is paid to obtaining new conditions for the monotonicity of the iterative process
for finding the lower coefficient of the parabolic equation, when the coefficient approaches from above.
This study continues the work~\cite{vabishchevich2017iterative}, where we consider iterative methods 
for the approximate solving the linear inverse problem of identifying the right-hand side 
for the parabolic equation.

The paper is organized as follows.
Statements of direct and inverse problems for the second-order parabolic equation are given in Section 2. 
The identification of the reaction coefficient that is independent of time is considered for the 
two-dimensional diffusion-reaction equation. An additional information on the solution of the equation 
is given at the final time moment.
An iterative adjustment algorithm for the desired coefficient is investigated in Section 3.
The proof of its monotonicity is based on the fulfillment of the maximum principle not only for the solution, 
but also for derivative of the solution with respect to time.
In Section 4, we construct a computational algorithm for approximate solving the identification problem 
for the lowest coefficient of the parabolic equation, and a discrete problem is formulated
using finite-element approximations in space and two-level time-stepping schemes.
Results of computational experiments for a model boundary-value problem are represented in Section 5.
The findings of the work are summarized in Section 6.

\section{Problem formulation}

The inverse problem of identifying the lowest coefficient of a parabolic equation is considered.
We confine ourselves to the two-dimensional case. Generalization to the 3D case is trivial. 
Let ${\bm x} = (x_1, x_2)$ and $\Omega$ be a bounded polygon.
The direct problem is formulated as follows.
We search $u({\bm x},t)$, $0 \leq t \leq T, \ T > 0$ such that 
it is the solution of the homogeneous parabolic equation of second order:
\begin{equation}\label{1}
  \frac{\partial u}{\partial t}- {\rm div} (k({\bm x}) \, {\rm grad} \, u) + c({\bm x}) u = f({\bm x},t),
  \quad {\bm x} \in \Omega,
  \quad 0 < t \leq T ,   
\end{equation} 
with coefficient $0 < k_1 \leq k({\bm x}) \leq k_2$.
The boundary conditions are also specified:
\begin{equation}\label{2}
  k({\bm x}) \frac{\partial u}{\partial n} + \mu({\bm x}) u = 0,
  \quad {\bm x} \in \partial\Omega,
  \quad 0 < t \leq T,    
\end{equation} 
where $\mu ({\bm x}) \geq 0, \  {\bm x} \in \partial \Omega$ and $\bm n$ is the normal to $\Omega$.
The initial conditions are
\begin{equation}\label{3}
  u({\bm x}, 0) = 0,
  \quad {\bm x} \in \Omega .  
\end{equation} 
The formulation (\ref{1})--(\ref{3}) presents the direct problem, 
where the coefficients of the equation as well as the boundary conditions are specified.

Let us consider the inverse problem, where in equation (\ref{1}), the lowest coefficient
$c({\bm x})$ that depends on spatial variables only is unknown. 
An additional condition is often formulated as 
\begin{equation}\label{4}
  u({\bm x}, T) = \psi ({\bm x}),
  \quad {\bm x} \in \Omega .  
\end{equation}
In this case, we have the case of the final redefinition.

Conditions for the unique solvability of the inverse coefficient problem (\ref{1})--(\ref{4})
and its correctness in various functional classes are established, for example,
in the works cited above (see~\cite{isakov1991inverse,prilepko1987solvability}).
We focus on using the iterative process to identify the coefficient
$c({\bm x})$, which has been employed, in particular, in~\cite{isakov2017inverse,kostin2015inverse,prilepko1992inverse}.
Let us formulate wider conditions for the monotonicity of the iterative process
of defining a new initial approximation, when the desired coefficient approaches from above.
In our consideration, we assume that the solution of the problem, the coefficients of the equation,
and the boundary conditions are sufficiently smooth, i.e., we have all necessary
derivatives with respect to the space variables and time.

On the set of functions satisfying the homogeneous boundary conditions (\ref{2}),
let us define the elliptic operator $\mathcal{A}$ by the relation
\[
 \mathcal{A} u = - {\rm div} (k({\bm x}) \, {\rm grad} \, u) .
\]
In this case, equations (\ref{1}), (\ref{2}) can be written in the compact form:
\begin{equation}\label{5}
  \frac{\partial u}{\partial t} + \mathcal{A} u + c({\bm x}) u = f({\bm x},t),
  \quad {\bm x} \in \Omega,
  \quad 0 < t \leq T .  
\end{equation}
Without loss of generality, we consider the inverse problem (\ref{3})--(\ref{5})
for the definition of the pair $(u,c)$ under a priori restrictions
on the reaction coefficient:
\begin{equation}\label{6}
  c({\bm x}) \geq 0,
 \quad \bm x \in \Omega . 
\end{equation} 
If $c({\bm x}) \geq m$ with a constant $m$, it is possible to employ the standard transition
to the problem for the function $v = \exp(m t) u$.

Assume that for the right-hand side of equation (\ref{1}) holds
\begin{equation}\label{7}
 f({\bm x}, 0) = 0,
 \quad \frac{\partial f}{\partial t}({\bm x}, t) >  0,  
 \quad \bm x \in \Omega .  
\end{equation}
Under these conditions, on the basis of the maximum principle for parabolic equations
(see, e.g., \cite{friedman2008partial,il1962linear}),
the solution $u$ at the final time moment is positive, i.e.
\begin{equation}\label{8}
 \psi ({\bm x}) > 0,
 \quad \bm x \in \Omega . 
\end{equation} 

\section{Iterative process}   

The inverse problem consists in evaluating the pair of functions $(u,c)$
from the conditions (\ref{3})--(\ref{5}) under the constraints (\ref{6})--(\ref{8}).
The iterative process of identifying the coefficient $c({\bm x})$ is implemented as follows.
It starts from specifying some initial approximation $c^0({\bm x})$.
With the known $c^k({\bm x})$, $k=0,1, ...$, where $k$ is the iteration number, the direct problem is solved:
\begin{equation}\label{9}
  \frac{\partial u^{k}}{\partial t} + \mathcal{A} u^k + c^k({\bm x}) u^k = f(\bm x,t),
  \quad {\bm x} \in \Omega,
  \quad 0 < t \leq T ,   
\end{equation} 
\begin{equation}\label{10}
  u^k({\bm x}, 0) = 0,
  \quad {\bm x} \in \Omega .  
\end{equation} 
A new approximation for the desired coefficient is evaluated from the equation at the final time moment $t=T$ 
using the redefinition (\ref{4}):
\begin{equation}\label{11}
  c^{k+1} ({\bm x}) \psi  = - \frac{\partial u^{k}}{\partial t}(\bm x,T) - \mathcal{A} \psi + f(\bm x,T),
  \quad {\bm x} \in \Omega .
\end{equation} 

In the works~\cite{isakov2017inverse,kostin2015inverse,prilepko1992inverse}, the initial approximation
is given in the form
\begin{equation}\label{12}
 c^0({\bm x}) = 0,
 \quad {\bm x} \in \Omega .  
\end{equation}
In this case, the monotone approach to the required coefficient ($c^{k+1}({\bm x}) \geq c^k({\bm x})$, 
approaching from below) holds, if this monotonicity condition holds for $k=1$: 
\begin{equation}\label{13}
 \mathcal{A} (u^{0} - \psi) \geq 0,
 \quad {\bm x} \in \Omega . 
\end{equation}
The condition (\ref{13}) is strong enough, but it can be removed.
To do this, consider the algorithm for monotone approaching the reaction coefficient $c({\bm x})$ from above.

For the initial-boundary value problem (\ref{3}), (\ref{5}), in assumption (\ref{7}), we have
\begin{equation}\label{14}
 u(\bm x,t) \geq 0,
 \quad \frac{\partial u}{\partial t} (\bm x,t) \geq 0,
 \quad {\bm x} \in \Omega,
 \quad 0 < t \leq T .  
\end{equation}
The non-negativity of the solution follows from the maximum principle and the non-negativity of 
the right-hand side ($f(\bm x,t) \geq 0$). The non-negativity of the time derivative is established 
similarly when considering
problem for ${\displaystyle w = \frac{\partial u}{\partial t}}$. Differentiation of  equation (\ref{5}) by time gives
\[
  \frac{\partial w}{\partial t} + \mathcal{A} w + c({\bm x}) w = \frac{\partial f}{\partial t} ({\bm x},t),
  \quad {\bm x} \in \Omega,
  \quad 0 < t \leq T .  
\]  
For $t=0$, from equation (\ref{5}) and the first condition in (\ref{7}), we get
\[
 w({\bm x}, 0) = 0,
 \quad \bm x \in \Omega .   
\]
From the maximum principle for this problem, it follows that $w({\bm x}, 0) \geq 0$.

In view of (\ref{14}), from equation (\ref{5}), for $t=T$, we obtain
\begin{equation}\label{15}
  c({\bm x}) \psi  \leq - \mathcal{A} \psi + f(\bm x,T),
  \quad {\bm x} \in \Omega .
\end{equation} 
Thus, the inverse problem (\ref{3})--(\ref{5}) is considered with two-side restrictions
(\ref{5}) and (\ref{15}) for the lowest coefficient $c({\bm x})$.

Let us consider the iterative process (\ref{9})--(\ref{11}) with the initial approximation
\begin{equation}\label{16}
  c^0({\bm x}) \psi  = - \mathcal{A} \psi + f(\bm x,T),
  \quad {\bm x} \in \Omega .
\end{equation} 
To find $u^0({\bm x},t)$, we solve the problem
\[
  \frac{\partial u^{0}}{\partial t} + \mathcal{A} u^0 + c^0({\bm x}) u^0 = f(\bm x,t),
  \quad {\bm x} \in \Omega,
  \quad 0 < t \leq T ,   
\]
\[
  u^0({\bm x}, 0) = 0,
  \quad {\bm x} \in \Omega .  
\] 
For $u^0({\bm x},t)$, similarly to (\ref{14}), we have
\[
 \frac{\partial u^0}{\partial t} (\bm x,t) \geq 0,
 \quad {\bm x} \in \Omega,
 \quad 0 < t \leq T .
\] 
In view of this, from (\ref{11}), we obtain
\[
  c^{1} ({\bm x}) \psi  = - \frac{\partial u^{0}}{\partial t}(\bm x,T) - \mathcal{A} \psi + f(\bm x,T)
 \leq - \mathcal{A} \psi + f(\bm x,T).
\] 
By (\ref{16}), we arrive at
\begin{equation}\label{17}
  c^{1} ({\bm x})  \leq  c^{0} ({\bm x}) ,
  \quad {\bm x} \in \Omega . 
\end{equation} 

Let us formulate a problem for the solution difference between two adjacent iterations:
\[
 \xi^{k}({\bm x}) = c^{k}({\bm x}) - c^{k-1}({\bm x}),
 \quad  w^{k}({\bm x},t) = u^{k}({\bm x},t) - u^{k-1}({\bm x},t) ,
 \quad k = 1,2, ... . 
\] 
From (\ref{9}), (\ref{10}), we have
\begin{equation}\label{18}
  \frac{\partial w^{k}}{\partial t} + \mathcal{A} w^k + c^{k-1} ({\bm x}) w^k = - \xi^k({\bm x}) u^k ,
\end{equation} 
\begin{equation}\label{19}
  w^k({\bm x}, 0) = 0.
\end{equation} 
From (\ref{11}), we get
\begin{equation}\label{20}
  \xi ^{k+1} ({\bm x}) \psi  = - \frac{\partial w^{k}}{\partial t} (\bm x,T).
\end{equation}

Similarly to (\ref{14}), we prove that
\begin{equation}\label{21}
 u^k(\bm x,t) \geq 0,
 \quad \frac{\partial u^k}{\partial t} (\bm x,t) \geq 0,
 \quad {\bm x} \in \Omega,
 \quad 0 < t \leq T ,  
\end{equation}  
for $k = 0,1, ...$. Considering the problem (\ref{18}), (\ref{19}) for $k=1$, in view of (\ref{17}),
on the basis of the maximum principle, we obtain
\begin{equation}\label{22}
 w^1 (\bm x,t) \geq 0,
 \quad \xi^1 (\bm x) \leq 0,
 \quad {\bm x} \in \Omega,
 \quad 0 < t \leq T . 
\end{equation}
An analogous property of the monotonicity of the approximate solution also holds for other $k = 2,3, ... $:
\begin{equation}\label{23}
 w^k (\bm x,t) \geq 0,
 \quad \xi^k (\bm x) \leq 0,
 \quad {\bm x} \in \Omega,
 \quad 0 < t \leq T .
\end{equation}  

The proof is by induction on $k$. For $k=1$, it is satisfied (see (\ref{22})).
Let us show that from the fulfillment (\ref{23}) for some $k$ this holds also for $k+1$.
If (\ref{23}) holds, taking into account the second inequality (\ref{21}),
after differentiating (\ref{18}) with respect to $t$, we obtain
\[
 \frac{\partial w^k}{\partial t} (\bm x,t) \geq 0,
 \quad {\bm x} \in \Omega,
 \quad 0 < t \leq T .  
\]
Under these conditions, directly from (\ref{20}), it follows that
\[
 \xi^{k+1} (\bm x) \leq  0,
 \quad {\bm x} \in \Omega ,  
\] 
and from (\ref{18}), (\ref{19}), for $k \rightarrow k+1$, we get
\[
 w^{k+1} (\bm x,t) \geq 0,
 \quad {\bm x} \in \Omega,
 \quad 0 < t \leq T .  
\] 

Define the error of the approximate solution as follows:
\[
 \delta c^{k}({\bm x}) = c^{k}({\bm x}) - c({\bm x}),
 \quad  \delta u^{k}({\bm x},t) = u^{k}({\bm x},t) - u({\bm x},t) ,
 \quad k = 0,1, ... . 
\] 
By (\ref{9}), (\ref{10}), we get
\begin{equation}\label{24}
  \frac{\partial \delta u^{k}}{\partial t} + \mathcal{A} \delta u^k + c ({\bm x}) \delta u^k = - \delta c^k({\bm x}) u ,
\end{equation} 
\begin{equation}\label{25}
  \delta u^k({\bm x}, 0) = 0.
\end{equation} 
From (\ref{11}), we have
\begin{equation}\label{26}
  \delta c^{k+1} ({\bm x}) \psi  = - \frac{\partial \delta u^{k}}{\partial t} (\bm x,T)
\end{equation}
for $k = 0,1, ... $.

For $k=0$, we obtain
\begin{equation}\label{27}
 \delta u^0 (\bm x,t) \leq  0,
 \quad \delta c^0 (\bm x) \geq  0,
 \quad {\bm x} \in \Omega,
 \quad 0 < t \leq T . 
\end{equation}
The second inequality follows immediately from (\ref{15}), (\ref{16}).
The first inequality is established for the solution of the problem (\ref{24}), (\ref{25}) with $k=0$
using the maximum principle. Further, similarly to (\ref{23}), on the basis of induction, 
the property of monotonicity is established for other $k = 1,2, ... $:
\begin{equation}\label{28}
 \delta u^k (\bm x,t) \leq  0,
 \quad \delta c^k (\bm x) \geq  0,
 \quad {\bm x} \in \Omega,
 \quad 0 < t \leq T .
\end{equation}  

The result of our consideration is the following statement on the monotonicity
of the iteration process (\ref{9})--(\ref{11}).

\begin{thm}\label{t-1}
The iteration process (\ref{9})--(\ref{11}) with the initial approximation specified by
(\ref{16}) is monotone and
\begin{equation}\label{29}
 \begin{split}
 u(\bm x,t) & \geq u^k (\bm x,t) \geq u^{k-1} (\bm x,t), \\
 c(\bm x) & \leq c^k (\bm x) \leq c^{k-1} (\bm x) , 
 \quad {\bm x} \in \Omega,
 \quad 0 < t \leq T ,
 \end{split}
\end{equation}  
for all $k = 1,2, ... $.
\end{thm}

If the initial condition is given in the form (\ref{12}), the similar statement that
\[
 \begin{split}
 u(\bm x,t) & \leq u^k (\bm x,t) \leq  u^{k-1} (\bm x,t), \\
 c(\bm x) & \geq c^k (\bm x) \geq  c^{k-1} (\bm x) ,
 \quad {\bm x} \in \Omega,
 \quad 0 < t \leq T , 
 \end{split}
\]
is proved under the additional condition (\ref{13}). 

\section{Computational implementation} 

It seems reasonable to recall some general points of
numerical solving the inverse coefficient problem (\ref{1})--(\ref{4}) 
on the basis of the iterative adjustment of the desired reaction coefficient.
The monotonicity of the iterative process is established in Theorem~\ref{t-1}
using the maximum principle for the solution and its time derivative (\ref{14}).
In constructing discretizations in space and time, we need to preserve this basic property of the
differential problem, i.e., an approximate solution of the problem should satisfy the maximum principle. 

Special attention should be given to monotone approximations in space 
(approximations of the diffusion-reaction operators) and discretizations in time.
The maximum principle is formulated in the most simple way (see, e.g., \cite{Samarskii}) 
for difference schemes on rectangular grids. For steady-state problems,
its implementation is associated with a diagonal dominance for the corresponding matrix and
non-positivity of off-diagonal elements. Some possibilities for constructing monotone approximations
on general irregular grids (using the finite volume method) and
the maximum principle for convection-reaction problems with anisotropic diffusion coefficients 
are discussed in the work~\cite{droniou2014finite}. 

Discretization in time leads to additional restrictions on monotonicity.
For example, a typical situation is the case, where the monotonicity of the approximate solution in two-level schemes 
is ensured by using a small enough step in time.
Unconditionally monotone time approximations for parabolic problems occure (see, for example, 
\cite{HundsdorferVerwer2003,Samarskii}) when using fully implicit two-level schemes (backward Euler).

Here, we focus on the application of the finite element method.
Monotone approximations in space for linear finite elements can be constructed
with restrictions on a computational grid (Delaunay-type mesh, see, e.g., \cite{huang2011discrete,letniowski1992three}).
Some additional restrictions arise (see, for instance, \cite{brandts2008discrete,ciarlet1973maximum})  
from  the reaction coefficient. They can be removed using the standard approach based on a correction 
of the approximations of the coefficient at the time derivative and reaction coefficient employing lumping procedures 
(see, e.g., \cite{chatzipantelidis2015preservation,Thomee2006}).

To solve numerically the problem (\ref{1})--(\ref{4}), we employ finite element 
approximations in space \cite{brenner2008mathematical,Thomee2006}. 
In the Hilbert space $H = L_2(\Omega)$, we define the scalar product and norm in the standard way:
\[
  (u,v) = \int_{\Omega} u({\bm x}) v({\bm x}) d{\bm x},
  \quad \|u\| = (u,u)^{1/2} .
\] 
We define the bilinear form
\[
 a(u,v) = \int_{\Omega } k \, {\rm grad} \, u \, {\rm grad} \, v  \, d {\bm x} +
 \int_{\partial \Omega } \mu \, u v d {\bm x} .
\] 
Define a subspace of finite elements $V^h \subset H^1(\Omega)$.
Let $\bm x_i, \ i = 1,2, ..., M_h$ be triangulation points for the domain $\Omega$.
When using Lagrange finite elements of the first order (piece-wise linear approximation),
we can define pyramid function $\chi_i(\bm x) \subset V^h, \ i = 1,2, ..., M_h$, where
\[
 \chi_i(\bm x_j) = \left \{
 \begin{array}{ll} 
 1, & \mathrm{if~}  i = j, \\
 0, & \mathrm{if~}  i \neq  j .
 \end{array}
 \right . 
\] 
For $v \in V_h$, we have
\[
 v(\bm x) = \sum_{i=i}^{M_h} v_i \chi_i(\bm x),
\] 
where $v_i = v(\bm x_i), \ i = 1,2, ..., M_h$.

Let us define a uniform grid in time  
\[
  t_n=n\tau,
  \quad n=0,1,...,N,
  \quad \tau N=T
\]
and denote $y_n = y(t_n), \ t_n = n \tau$.
Define an approximate solution of the inverse problem (\ref{1})--(\ref{4}) as
\[
 w_{n}(\bm x) \in V^h, \quad n=0,1,...,N,
 \quad s(\bm x) \in V^h . 
\]
For the fully implicit scheme, the solution is evaluated from
\begin{equation}\label{30}
 \begin{split}
 & \left (\frac{w_{n+1} - w_{n}}{\tau} , v \right ) + a(w_{n+1}, v) + (s w_{n+1}, v) = (f_{n+1}, v), \\
 & \qquad n=0,1,...,N-1,
 \quad 0 < t \leq T, 
 \quad v \in V^h . 
 \end{split}
\end{equation} 
From (\ref{3}), (\ref{4}), we have
\begin{equation}\label{31}
 w_{0} = 0,
 \quad (w_{N}, v) = (\psi, v),
 \quad v \in V^h .   
\end{equation} 

The computational algorithm for solving the problem (\ref{30}), (\ref{31}) is based on the iterative
method (\ref{9})--(\ref{11}), (\ref{16}) for identifying the lowest coefficient.
The calculation starts from specifying the initial approximation for the desired coefficient:
\begin{equation}\label{32}
  (c^0({\bm x}) \psi, v) = -  a(\psi, v)  + (f(\bm x,T), v),
 \quad v \in V^h .  
\end{equation}  
For known $c^k({\bm x})$, we solve the direct problem for evaluating $w_{n}^k({\bm x})$:
\begin{equation}\label{33}
 \begin{split}
 & \left (\frac{w_{n+1}^k - w_{n}^k}{\tau} , v \right ) + a(w_{n+1}^k, v) + (s^k w_{n+1}^k, v) = (f_{n+1}, v), \\
 & \qquad n=0,1,...,N-1,
 \quad 0 < t \leq T, 
 \end{split}
\end{equation} 
\begin{equation}\label{34}
 w_{0}^k = 0,
 \quad v \in V^h .   
\end{equation}
After this, the reaction coefficient is adjusted:
\begin{equation}\label{35}
  (c^{k+1}({\bm x}) \psi, v) = - \frac{w_{N}^k - w_{N-1}^k}{\tau } -  a(\psi, v)  + (f(\bm x,T), v),
 \quad v \in V^h .  
\end{equation} 

If we apply (\ref{32})--(\ref{35}), no additional procedures are needed for monotonization.

\section{Numerical experiments} 

To illustrate the capabilities of the iterative technique for solving inverse problems of identifying 
the lowest coefficient of parabolic equations, we present the results of numerical experiments for a test problem.
Let us consider model problem (\ref{1})--(\ref{3}), where
\[
 k({\bm x}) = 1,
 \quad f({\bm x}, t) = 100 t \exp(-x_1),
 \quad \mu({\bm x}) = 10,
 \quad T = 0.25 . 
\] 
The problem is solved in the unit square
\[
 \Omega = \{ \bm x = (x_1, x_2) \ | \ 0 < x_1 < 1, \ 0 < x_2 < 1 \} 
\]
The data at the final time moment (see (\ref{4})) are obtained
from the solution of the direct problem with a given coefficient $c({\bm x})$.

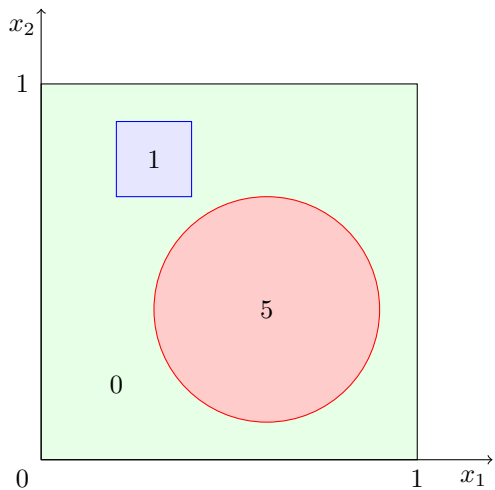
\begin{figure}[ht] 
  \begin{center}
    \begin{tikzpicture}[scale=0.5]
	\fill[color=green!10,draw=black] (0,0) rectangle +(10,10);
	\fill[color=red!20, draw=red] (6,4) circle (3);
	\fill[color=blue!10,draw=blue] (2,7) rectangle +(2,2);
        \draw [->] (0,0) -- (12,0);
        \draw [->] (0,0) -- (0,12);
        \draw(-0.5,-0.5) node {$0$};
        \draw(10,0-0.5) node {$1$};
        \draw(11.5,-0.5) node {$x_1$};
        \draw(-0.5,10) node {$1$};
        \draw(-0.5,11.5) node {$x_2$};
        \draw(6,4) node {$5$};
        \draw(3,8) node {$1$};
        \draw(2,2) node {$0$};
   \end{tikzpicture}
    \caption{Reaction coefficient}
   \label{f-1}
  \end{center}
\end{figure}

In our case, the coefficient $c({\bm x})$ is the piecewise constant (see Fig.~\ref{f-1}):
inside a circle of radius 0.3 with the center (0.6,0.4), we put $c({\bm x})=5$;
inside the square with side 0.2 and the center (0.3,0.8), we have $c({\bm x})=1$;
and otherwise, we put $c({\bm x})=0$.

The solution of the direct problem with this coefficient $c({\bm x})$ at the final time moment
(the function $u(\bm x,T)$) is used as input data for the inverse problem. 
In our analysis, we focus on iterative solving the identification problem after finite element discretizations in space.
Because of this, we do not discuss the dependence of the accuracy of the numerical solution on approximations in space,
it seems appropriate to do in a separate study. The effect of computational errors is studied via 
calculations on different time grids, when the input data is derived from the solution of the direct problem 
on more fine time grids and with higher-order approximations in time.

\begin{figure}[htp]
  \begin{center}
    \includegraphics[scale = 0.75] {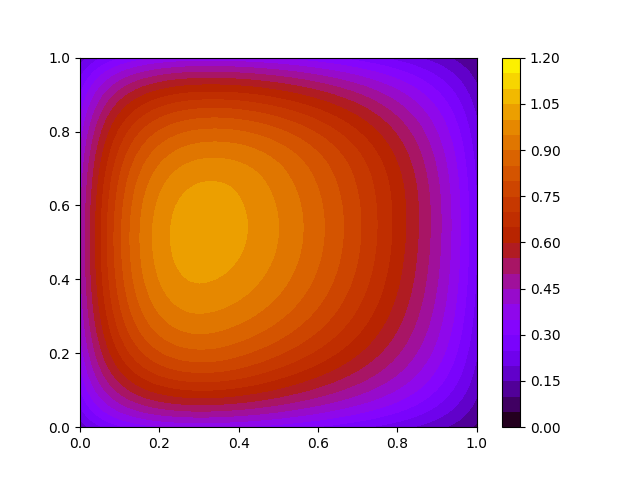}
	\caption{The solution of the direct problem $u(\bm x,T)$:
         $u_{\min} = 0.0884557$,  $u_{\max} = 1.03433$}
	\label{f-2}
  \end{center}
\end{figure}

To solve the direct problem, we employ the time step $\tau = 1\cdot 10^{-5}$.
The division into 50 intervals in each direction is used to construct the uniform spatial grid, 
the Lagrangian finite elements of first degree are applied.
The solution at the finite time moment is shown in Fig.~\ref{f-2}.

The inverse problem is solved using the fully implicit scheme (see (\ref{33})--(\ref{35})).
The error of the approximate solution of the identification problem on a separate iteration is evaluated as follows:
\[
 \varepsilon_\infty(k) = \max_{\bm x \in \Omega} |c^k(\bm x) - c(\bm x)|,
\] 
\[
 \varepsilon_2(k) = \|c^k(\bm x) - c(\bm x)\| .
\] 

\begin{figure}[htp]
  \begin{center}
\begin{minipage}{0.49\linewidth}
\center{\includegraphics[width=1\linewidth]{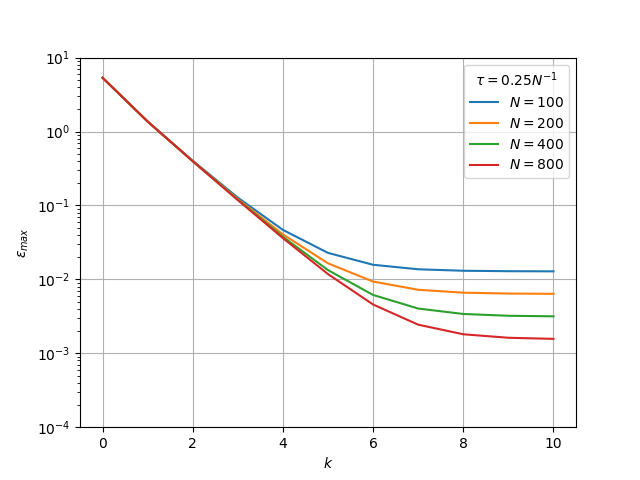}} \\
\end{minipage}
\hfill
\begin{minipage}{0.49\linewidth}
\center{\includegraphics[width=1\linewidth]{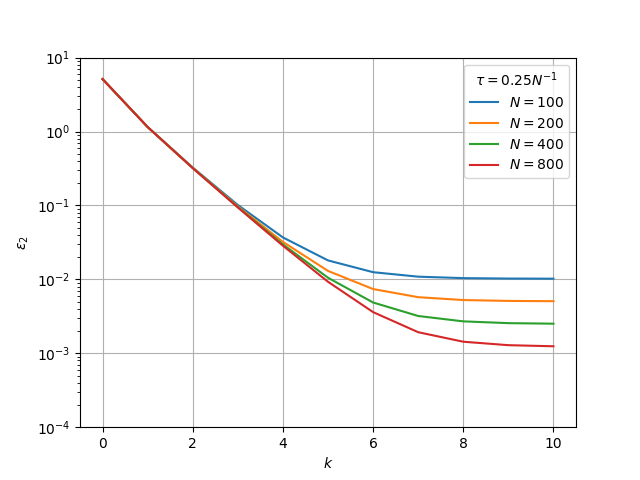}} \\
\end{minipage}
\caption{The iterative process with the initial approximation (\ref{32}) } 
\label{f-3}
  \end{center}
\end{figure}

\begin{figure}[htp]
  \begin{center}
\begin{minipage}{0.49\linewidth}
\center{\includegraphics[width=1\linewidth]{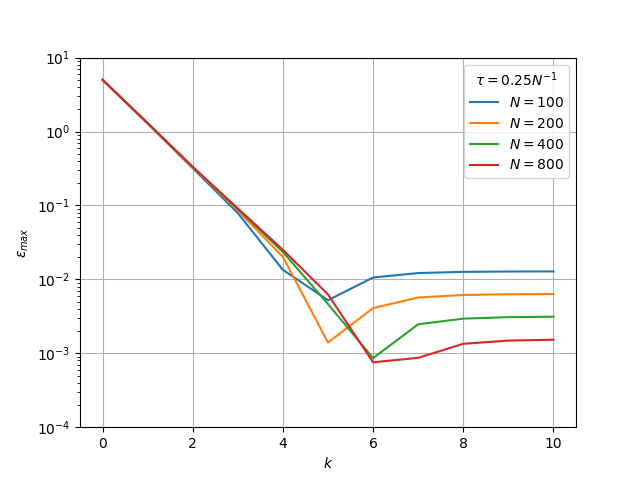}} \\
\end{minipage}
\hfill
\begin{minipage}{0.49\linewidth}
\center{\includegraphics[width=1\linewidth]{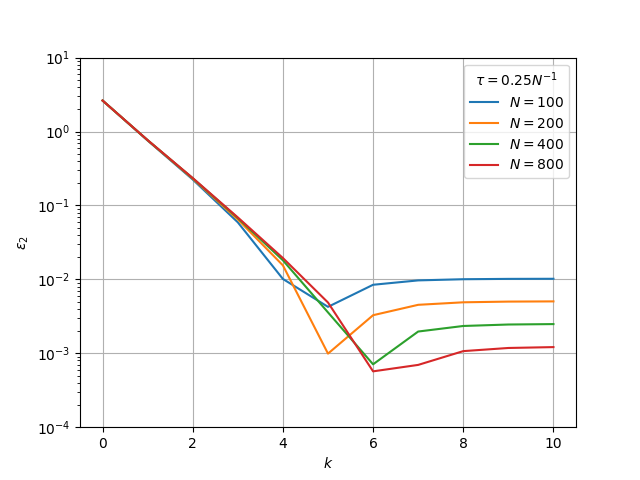}} \\
\end{minipage}
\caption{The iterative process with the initial approximation (\ref{12})} 
\label{f-4}
  \end{center}
\end{figure}

The main issue is to evaluate an actual convergence rate for the iterative processes under the consideration. 
We need to recognize clearly how quickly the accuracy of the approximate solution is stabilized with increasing
the iteration number. The obtained error itself depends on a time step, namely,
the smaller time step, the higher accuracy of the approximate solution.
Influence of the time step of the iterative process (\ref{33})--(\ref{35}) 
with the initial approximation (\ref{32}) on accuracy is shown in Fig.~\ref{f-3}.
We observe a high convergence of the iterative process and the improvement 
of the accuracy of the approximate solution by reducing the time step. 
Similar results for the iterative process with the initial approximation
(\ref{12}) are presented in Fig.~\ref{f-4}. 

\begin{figure}[htp]
\begin{minipage}{0.33\linewidth}
  \begin{center}
    \includegraphics[width=1.\linewidth] {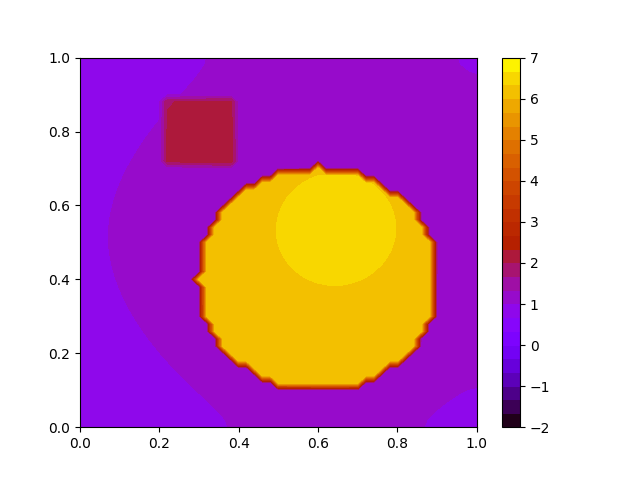}
  \end{center}
\end{minipage}\hfill
\begin{minipage}{0.33\linewidth}
  \begin{center}
    \includegraphics[width=1.\linewidth] {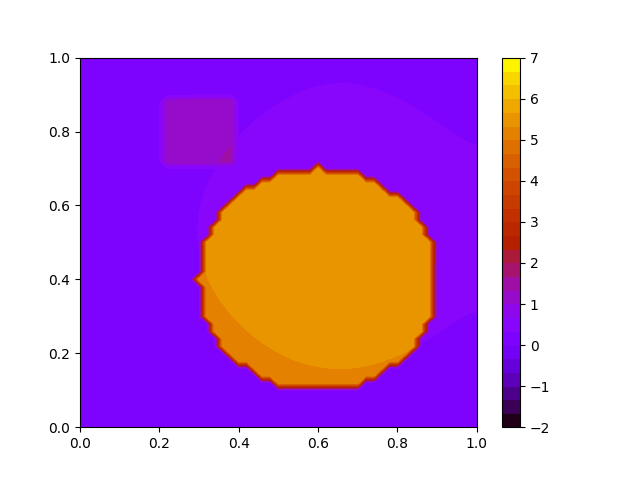}
  \end{center}
\end{minipage}
\begin{minipage}{0.33\linewidth}
  \begin{center}
    \includegraphics[width=1.\linewidth] {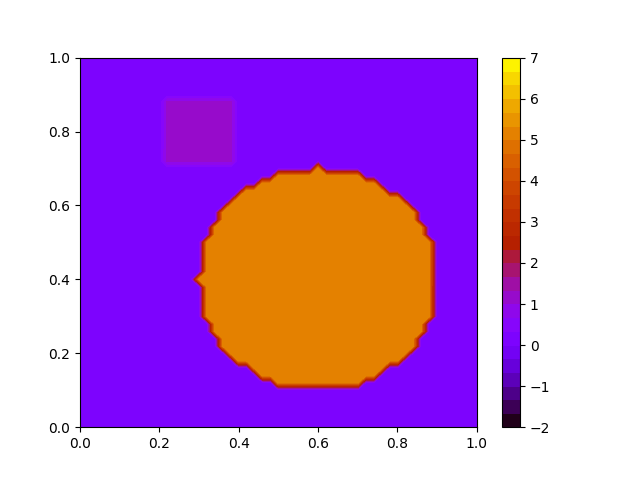}
  \end{center}
\end{minipage}
	\caption{The solution $c({\bm x})$  for different $k$ with initial approximation (\ref{32}): left --- $k = 1$, center ---  $k = 2$,
         right --- $k = 3$}
	\label{f-5}
\end{figure}

\begin{figure}[htp]
\begin{minipage}{0.33\linewidth}
  \begin{center}
    \includegraphics[width=1.\linewidth] {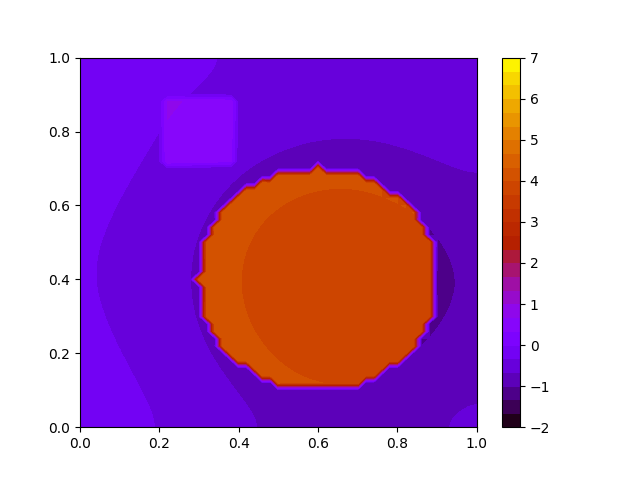}
  \end{center}
\end{minipage}\hfill
\begin{minipage}{0.33\linewidth}
  \begin{center}
    \includegraphics[width=1.\linewidth] {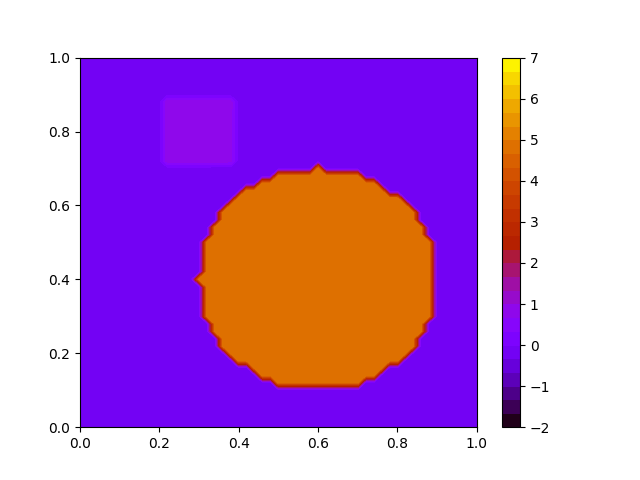}
  \end{center}
\end{minipage}
\begin{minipage}{0.33\linewidth}
  \begin{center}
    \includegraphics[width=1.\linewidth] {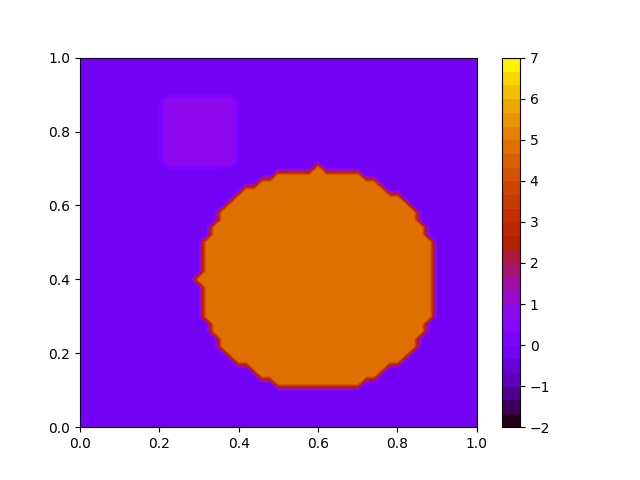}
  \end{center}
\end{minipage}
	\caption{The solution $c({\bm x})$  for different $k$ with initial approximation (\ref{12}): left --- $k = 1$, center ---  $k = 2$,
         right --- $k = 3$}
	\label{f-6}
\end{figure}

The convergence of the approximate solution for the reaction coefficient is shown in Fig.~\ref{f-5} and Fig.~\ref{f-6}. 
For the iterative process with initial approximation (\ref{32}),
we observe a monotone convergence from above (see Fig.~\ref{f-5}).
The iterative process with the initial approximation (\ref{12}) is non-monotone.
In particular, for $k=1$, on a part of the domain $\Omega$,
the function $c^1({\bm x})$ is negative (see Fig.~\ref{f-6}).

\section{Conclusion} 

\begin{enumerate}
 \item A nonlinear inverse problem of identifying the lowest coefficient that depends only on spatial variables
  is studied for a second-order parabolic equation. The solution of the parabolic equation at the final time moment 
  is given, i.e., the final redefinition is considered.
 \item An iterative process of identifying an unknown coefficient is conducted by solving the standard
  initial-boundary value problem at each iteration. The main result is in establishing the
  monotonicity of the iterative process, where the desired lower coefficient approaches from above.
 \item The computational algorithm is based on standard approximations in space by linear finite elements,
  whereas time-stepping is implemented using the fully implicit two-level schemes.
 \item Possibilities of the proposed algorithms were demonstrated by numerical solving a test 
  two-dimensional problem.
\end{enumerate} 

\section*{Acknowledgements}

The work was supported by the mega-grant of the Russian Federation Government 
(\# 14.Y26.31.0013) and RFBR (project 17-01-00689).


\end{document}